# Explicit Upper Bound for the Rank of $J_0(q)$

E. Kowalski, Princeton University   P. Michel, Université d'Orsay

**Abstract**

We refine the techniques of our previous papers [KM1], [KM2] to prove that the average order of vanishing of $L$-functions of primitive automorphic forms of weight 2 and prime level $q$ satisfies

$$\sum_{f \in S_2(q)^*} \mathrm{ord}_{s=\frac{1}{2}} L(f,s) \leq C |S_2(q)^*|$$

with $C = 6.5$, for all $q$ large enough. On the Birch and Swinnerton-Dyer conjecture, this implies

$$\mathrm{rank}\, J_0(q) \leq C \dim J_0(q)$$

for $q$ prime large enough.

## Contents

## 1  Introduction

In the papers [KM1] and [KM2], we have shown, among other results, that the analytic rank of the abelian variety $J_0(q) = \mathrm{Jac}\, X_0(q)$ (i.e. the order of vanishing at the central critical point of their Hasse-Weil $L$-functions, see [MSD]) satisfies

$$\mathrm{rank}_a J_0(q) \leq C \dim J_0(q)$$

for all prime numbers $q$, $C$ being some absolute constant. It is our purpose here to show how to compute an admissible value of $C$. To put the result in perspective, recall that it is conjectured [Br1] that

$$\mathrm{rank} J_0(q) = \mathrm{rank}_a J_0(q) \sim \frac{1}{2} \dim J_0(q)$$

(based on the consideration of the sign of the functional equation of automorphic $L$-functions), the first equality being the famous conjecture of Birch and Swinnerton-Dyer in this particular case. Assuming the Riemann Hypothesis for automorphic $L$-functions, Iwaniec, Luo and Sarnak have recently proved that one could take $C = \frac{99}{100}$; the best known previously was $C = \frac{23}{22}$ [KM1], or $C = 1$ (using also the Riemann Hypothesis for Dirichlet $L$-functions).

**Theorem 1** *Assume that the Birch and Swinnerton-Dyer conjecture holds for the abelian varieties $J_0(q)$, $q$ prime. Then for $q$ large enough we have*

$$\mathrm{rank} J_0(q) \leq 6.5 \dim J_0(q). \tag{1}$$

By Eichler-Shimura theory, we have a factorization

$$L(J_0(q), s) = \prod_{f \in S_2(q)^*} L(f, s)$$

where the product is over the finite set (of cardinality $|S_2(q)^*| = \dim J_0(q)$) of all primitive weight 2 forms of level $q$. Hence, assuming the Birch and Swinnerton-Dyer conjecture for the $J_0(q)$, Theorem 1 is equivalent with

**Theorem 2** *For any prime number $q$ large enough, we have*

$$\sum_{f \in S_2(q)^*} \mathrm{ord}_{s=\frac{1}{2}} L(f, s) \leq 6.5 |S_2(q)^*|$$

*(here $s = \frac{1}{2}$ is the central critical point in the analytic normalization of automorphic $L$-functions).*

The proof essentially follows that of [KM1], [KM2], keeping track of the constant involved. To get a good bound, we refine somewhat the previous method, and this requires some very technical arguments. We will present the main steps in outline, leading to the computation of $C$; all technical results are reserved for the later sections of the paper, or for the appendices. We emphasize that at some points we have made (almost) arbitrary choices of some of the parameters involved. It is likely that the constant can be slightly improved by other adjustments.

**Acknowledgments.** This article has benefited from encouragements and suggestions from É. Fouvry, H. Iwaniec and P. Sarnak. The computations were performed on the PARI System of Batut, Bernardi, Cohen.

**Conventions.** In the following, $\varepsilon$ will usually be used to denote a real number $> 0$ with the understanding that for $\varepsilon > 0$ the stated inequality holds. Similarly for $A$, $B$, which are understood to be real numbers $> 0$ which are stated to exist such that the inequality holds. We note $\log_n$ the $n$-th iterate of log.

The function $E(x)$ is defined by

$$E(x) = \int_1^{+\infty} \frac{e^{-tx}}{t^2} dt = x \int_x^{+\infty} y^{-1} e^{-y} \frac{dy}{y} = x\Gamma(-1, x), \qquad (2)$$

(in terms of the incomplete Gamma function).

## 2 Computation of the constant

We recall the necessary notations. For any $f \in S_2(q)^*$, its Fourier expansion is written as

$$f(z) = \sum_{n \geq 1} \lambda_f(n) n^{1/2} e(nz)$$

(where, as usual, $e(z) = \exp(2\pi i z)$), and its Hecke $L$-function is

$$L(f, s) = \sum_{n \geq 1} \lambda_f(n) n^{-s}$$

with $\lambda_f(1) = 1$. The completed $L$-function is

$$\Lambda(f, s) = \left(\frac{\sqrt{q}}{2\pi}\right)^s \Gamma(s + \tfrac{1}{2}) L(f, s)$$

which satisfies the functional equation

$$\Lambda(f, s) = \varepsilon_f \Lambda(f, 1-s)$$

with $\varepsilon_f = -q^{1/2} \lambda_f(q)$.

The Petersson inner-product is denoted by $(f, g)$ for two forms $f$ and $g$. We use the symbol $\sum_f^h$ to denote summations over modular forms with the so-called harmonic weight $\omega_f = \frac{1}{4\pi(f,f)}$ inserted, i.e.

$$\sum_{f \in S_2(q)^*}^h \alpha_f := \sum_{f \in S_2(q)^*} \frac{\alpha_f}{4\pi(f, f)}.$$

### 2.1 Step 1: harmonic average

Instead of Theorem 2, it is enough to prove

**Theorem 3** *For $q$ prime, large enough, we have*

$$\sum_{f \in S_2(q)^*}^h \mathrm{ord}_{s=\frac{1}{2}} L(f, s) \leq 6.5,$$

Indeed, the corresponding unweighted inequality

$$\sum_{f \in S_2(q)^*} \mathrm{ord}_{s=\frac{1}{2}} L(f, s) \leq 6.5 \dim J_0(q)$$

can be obtained by the same technique for removing the harmonic weight already employed in [KM2, Section 5], [KM3] (see also [Kow, Ch. 3] for a somewhat more detailed treatment).

## 2.2 Step 2: the explicit formula

Let $\varphi$ be a test function, $C^\infty$, even, with support in $[-1, 1]$, such that $\varphi(0) = 1$ and which is such that the Laplace transform $\hat{\varphi}$ of $\varphi$, defined by

$$\hat{\varphi}(s) = \int_{\mathbf{R}} f(x)e^{sx}dx$$

(which is an entire function) satisfies the positivity condition

$$\hat{\varphi}(0) = \int_{\mathbf{R}} f(x)dx > 0, \text{ and } \operatorname{Re} \hat{\varphi}(s) \geq 0 \text{ for } |\operatorname{Re}(s)| \leq 1. \tag{3}$$

Let $\lambda$ be a parameter, $\lambda = \theta \log q$ for some parameter $\theta > 0$ which will be fixed later on. We let

$$\varphi_\lambda(x) = \varphi\left(\frac{x}{\lambda}\right)$$

so that

$$\hat{\varphi}_\lambda(s) = \lambda \hat{\varphi}(\lambda s).$$

By successive integrations by parts we have the inequalities, for all integer $k \geq 0$ and all $s \in \mathbf{C}$:

$$|\hat{\varphi}_\lambda(s)| \leq ||\varphi^{(k)}||_1 \frac{\lambda}{|\lambda s|^k} e^{\lambda |\operatorname{Re}(s)|} \tag{4}$$

where $|| \cdot ||_1$ is the $L^1$-norm.

**Proposition 1** *Let $\theta$ be such that $0 < \theta < \frac{3}{2}$. Then we have*

$$\sum_{f \in S_2(q)^*}^h \operatorname{ord}_{s=\frac{1}{2}} L(f, s) \leq \frac{1}{\theta \hat{\varphi}(0)} + \frac{1}{2} + \frac{2}{\lambda \hat{\varphi}(0)} Z + o(1) \tag{5}$$

*where*

$$Z = \sum_{f \in S_2(q)^*}^h \sum_{\substack{L(f,\rho)=0 \\ \beta - \frac{1}{2} > \lambda^{-1}}} \left| \hat{\varphi}_\lambda(\rho - \tfrac{1}{2}) \right|. \tag{6}$$

This is the consequence of the explicit formulae of Riemann-Weil-Mestre, already used in this context by Brumer and Murty. See Section 3.1 for a short account.

Let $T > 0$ be another parameter to will be chosen later. We write the sum $Z$ as $Z = Z_T + Z^T$, where $Z_T$ is that part where $\rho$ runs over the zeros with $|\gamma| \leq T$ and $Z^T$ is the part where $|\gamma| > T$. We further write $Z_T = Z'_T + Z''_T$, where in $Z''_T$ we consider those zeros $\rho = \beta + i\gamma$ with $\beta - \frac{1}{2} > \frac{\log_2 q}{\log q}$, and therefore $Z'_T$ is the sum containing the (hypothetical) zeros very near to $\frac{1}{2}$. As the intuition suggests, $Z'_T$ will be the most difficult term to manage.

## 2.3 Step 3: zeros "far" from the critical point

**Proposition 2** *Assume that $T = \log_2 q / \log q$. Then we have*

$$Z^T = o(\lambda), \quad Z'_T = o(\lambda).$$

See Section 3.2 for the proof.

## 2.4 Step 4: zeros close to the critical point

Let $\lambda = a\Delta \log q$, and subdivide the rectangle $[\frac{1}{2} + \frac{1}{\lambda}, \frac{1}{2} + \frac{\log_2 q}{\log q}] \times [-T, T]$ into $O(\log_2 q)$ rectangles of the form

$$[\tfrac{1}{2} + \tfrac{1}{\lambda}, \tfrac{1}{2} + \tfrac{\log_2 q}{\log q}] \times [\tfrac{n}{\lambda}, \tfrac{n+1}{\lambda}]$$

with $|n| \ll \log_2 q$ ($T$ is now assumed to be chosen as in Proposition 2).

For any fixed $n \geq 0$, let

$$Z_T(n) = \sum_{f \in S_2(q)^*}^{h} \sum_{\substack{\beta - 1/2 > \lambda^{-1} \\ n/\lambda \leq \gamma < (n+1)/\lambda}} \left|\hat{\varphi}_\lambda(\rho - \tfrac{1}{2})\right|$$

so

$$Z'_T \leq \sum_n Z_T(n). \tag{7}$$

**Proposition 3** *Let $0 < \Delta < \frac{1}{2}$, $M = q^\Delta$, $\lambda = a\Delta \log q$ with $0 < a < 1$, so $\theta = a\Delta$. Let $\psi$ be any of the three functions*

$$|\varphi|, \quad \frac{1}{n}|\varphi'|, \quad \frac{1}{n^2}|\varphi''|,$$

*let $c = 4\pi \sin \frac{\pi-1}{2} = 11.028...$, let $F(a, u)$ be the function (see (2) for $E$)*

$$F(a, u) = \frac{1}{cu^2}\Big(e^{-2u/a} + ue^{-u}E(\tfrac{2-a}{a}u) - ue^u E(\tfrac{2+a}{a}u)\Big),$$

*let $K(a, x)$ be the function*

$$K(a,x) = \frac{2}{c}\Big\{E\big(\tfrac{1}{2}(\tfrac{2}{a} - x)\big) + \frac{1}{x-1}\Big(E\big(\tfrac{1}{2}(\tfrac{2}{a} - x)\big) - e^{(x-1)/2}E\big(\tfrac{1}{2}(\tfrac{2}{a} - 1)\big)\Big)$$
$$- \frac{1}{x+1}\Big(E\big(\tfrac{1}{2}(\tfrac{2}{a} - x)\big) - e^{(x+1)/2}E\big(\tfrac{1}{2}(\tfrac{2}{a} + 1)\big)\Big)\Big\},$$

*let $G_\psi(a)$ be defined by*

$$G_\psi(a) = F(a, \tfrac{1}{2})\hat{\psi}(1) + \int_{-1}^{1} x\psi(x)K(a,x)e^{x/2}dx.$$

*Then we have, as $q \to +\infty$*

$$\frac{1}{\lambda}Z_T(n) \leq \frac{a^2}{(1-a)^2}(G_\psi(1) - G_\psi(a)) + o_{a,\varphi}(\tfrac{1}{n^j}). \tag{8}$$

For the proof, see Section 3.3.

## 2.5 Step 5: choosing the test function

For any $\varepsilon > 0$, we let $g_\varepsilon$ be any even real-valued smooth function on $\mathbf{R}$, which is strictly increasing for $x < 0$ and satisfies

$$g_\varepsilon(x) = 0, \text{ for } |x| > \tfrac{1}{2} + \varepsilon, \quad g_\varepsilon(x) = 1, \text{ for } |x| \leq \tfrac{1}{2}$$

and we let

$$\varphi_\varepsilon = \frac{g_\varepsilon \star g_\varepsilon}{\cosh}.$$

Then $\varphi_\varepsilon$ satisfies all the assumptions described at the beginning of Section 2.2. Moreover, for $\varepsilon \to 0$, it is clear that $\varphi_\varepsilon$ converges, in the sense of distribution, to $\varphi_0$, where

$$\varphi_0(x) = \frac{\min\{0, 1 - |x|\}}{\cosh(x)}.$$

**Proposition 4** *With notations and parameters chosen as in Propositions 2 and 3, we have as $q \to +\infty$*

$$\frac{1}{\lambda} Z'_T \leq 2 \frac{a^2}{(1-a)^2} \lim_{\varepsilon \to 0} \left\{ 3(G_{\varphi_\varepsilon}(1) - G_{\varphi_\varepsilon}(a)) + \left(\frac{\pi^2}{6} - \frac{5}{4}\right)(G_{\varphi''_\varepsilon}(1) - G_{\varphi''_\varepsilon}(a)) \right\} + o_a(1).$$

Now we recapitulate the computations leading to our value of $C$:

1. For $\varphi_0$, we have

$$\hat{\varphi}_0(0) = 2 \int_0^1 \frac{1-x}{\cosh x} dx = 0.9281...$$

2. By Proposition 1, Proposition 2, Proposition 3 and Proposition 4

$$\text{rank}_a J_0(q) \leq H(a) + \eta + o(1) \tag{9}$$

for any $\eta > 0$ where $H(a)$ is equal to

$$\frac{1}{2} + \frac{1}{\hat{\varphi}(0)} \left( \frac{2}{a} + \frac{4a^2}{(1-a)^2} \lim_{\varepsilon \to 0} 3(G_{\varphi_\varepsilon}(1) - G_{\varphi_\varepsilon}(a)) + (\pi^2/6 - 5/4)(G_{\varphi''_\varepsilon}(1) - G_{\varphi''_\varepsilon}(a)) \right)$$

and $a$, $0 < a < 1$, can be chosen at will (we have taken $\Delta = \tfrac{1}{2}$, which by continuity explains the $\eta$ occuring on the right-hand side).

3. By Lemma 3 below, we have

$$G_{|\varphi_\varepsilon|}(1) \to 0.1535...$$
$$G_{|\varphi''_\varepsilon|}(1) \to 0.3321...$$

as $\varepsilon \to 0$.

4. Similarly, we can compute, for each value of $a$, the limits

$$\lim_{\varepsilon \to 0^+} G_{|\varphi_\varepsilon|}(a), \quad \lim_{\varepsilon \to 0^+} G_{|\varphi''_\varepsilon|}(a)$$

using Lemma 2 and the definition in Proposition 3.

5. We choose $a = 0.48$ (after a few numerical experimentations with (9), with no attempt to optimize beyond the second decimal). We get, for some $\eta > 0$ small enough, by computing $H(0.48) \simeq 6.498$
$$H(0.48) + \eta \leq 6.5.$$

Thus Theorem 1 is proved.

**Remark.** A file containing the (commented) PARI/GP programs used for the computations described above is available from the authors upon request.

## 3 The intermediate statements

### 3.1 Proposition 1

By reasoning as in [KM1], we have the inequality (for $q \to +\infty$)

$$\lambda \hat{\varphi}(0) \sum_{f \in S_2(q)^*}^{h} \mathrm{ord}_{s=\frac{1}{2}} L(f,s) \leq (1+o(1)) \log q - 2 \sum_{f \in S_2(q)^*}^{h} S_{1,f} - 2 \sum_{f \in S_2(q)^*}^{h} S_{2,f}$$
$$+ 2 \sum_{f \in S_2(q)^*}^{h} \sum_{\substack{L(f,\rho)=0 \\ \beta - \frac{1}{2} > \lambda^{-1}}} \left| \hat{\varphi}_\lambda(\rho - \tfrac{1}{2}) \right|$$

with

$$S_{1,f} = \sum_p \lambda_f(p) \frac{\varphi_\lambda(\log p) \log p}{\sqrt{p}}$$
$$S_{2,f} = \sum_p (\lambda_f(p^2) - 2) \frac{\varphi_\lambda(2 \log p) \log p}{\sqrt{p}}.$$

Standard estimates yield

$$\sum_{f \in S_2(q)^*}^{h} S_{1,f} = O(e^\lambda q^{-3/2}) + O(1)$$

and

$$\sum_{f \in S_2(q)^*}^{h} S_{2,f} = \frac{\lambda}{2} \hat{\varphi}(0) + O(e^{\lambda/2} q^{-3/2}) + O(1).$$

Hence the result.

### 3.2 Proposition 2

As usual, we let $N(f; \sigma, t_1, t_2)$ denote the number of non-trivial zeros of $L(f,s)$ contained in the rectangle $[\sigma, 1] \times [t_1, t_2]$.

To estimate $Z^T$ and $Z'_T$, we use the following density theorem for the zeros of automorphic $L$-functions, a strengthened version of Theorem 1.3 of [KM2], which is the case where $\Delta < 1/4$.

**Theorem 4** *There exists an absolute constant $B > 0$ such that for any $\Delta < \frac{1}{2}$, any $\sigma \geq \frac{1}{2} + (\log q)^{-1}$, and any $t_1 < t_2$ with $t_2 - t_1 \geq (\log q)^{-1}$, we have*

$$\sum_{f \in S_2(q)^*}^{h} N(f, \sigma, t_1, t_2) \ll_\Delta (1 + |t_1| + |t_2|)^B q^{-2\Delta(\sigma - \frac{1}{2})} (t_2 - t_1) \log q.$$

For the proof of this result, see Appendix A

Recall that

$$Z^T = \sum_{f \in S_2(q)^*}^{h} \sum_{\substack{L(f,\rho)=0 \\ \beta - \frac{1}{2} > \lambda^{-1}, |\gamma| > T}} \left| \hat{\varphi}_\lambda(\rho - \tfrac{1}{2}) \right|$$

and similarly $Z_T''$ is the sum over zeros $\rho = \beta + i\gamma$ with $\beta - \frac{1}{2} > \frac{\log_2 q}{\log q}$, $|\gamma| \leq T$.

We subdivide the two strips defined by

$$z = x + iy, \quad \frac{1}{2} + \frac{1}{\lambda} < x \leq 1, \quad |y| \geq T$$

into small squares of size $\lambda^{-1}$. From Theorem 4 and from (4) we obtain

$$\frac{1}{\lambda} Z^T \ll_{k,\Delta} \sum_{1 \leq m \leq \lambda} \sum_{n \geq 0} \frac{(1 + T + n/\lambda)^B}{(\lambda T + n)^k} e^{m(1-c/\theta)} \ll_{k,\Delta,\theta} (\lambda T)^{-(k-1-B)}$$

if $\theta < \Delta < \frac{1}{2}$. For $T = \log_2 q / \log q$ and then $k$ large enough so that the resulting bound gives

$$Z^T = o(\lambda) \tag{10}$$

as $q$ tends to infinity.

Similarly, we see that

$$Z_T'' = o(\lambda). \tag{11}$$

### 3.3 Proposition 3

In Appendix B, we will establish another density theorem, more precise than Theorem 4 close to the critical line.

**Theorem 5** *Let notations be as in the statement of Proposition 3. Then for any $\sigma$ such that*

$$\frac{1}{a \log M} \leq \sigma - \frac{1}{2} \leq \frac{\log_2 q}{\log q}$$

*and any $t_1$, $t_2$ such that*

$$-\frac{\log_2 q}{\log q} \leq t_1 < t_2 \leq \frac{\log_2 q}{\log q}$$

*and $t_2 - t_1 = (a \log M)^{-1}$, we have*

$$\sum_{f \in S_2(q)^*}^{h} N(f; \sigma, t_1, t_2) \leq \frac{a^2(1+o(1))}{(1-a)^2} (F(1,u) - F(a,u)) + O_{a,\Delta}(\exp(-2u) \frac{(\log_2 q)^3}{\log q}).$$

*where we have written $u = \lambda(\sigma - \frac{1}{2}) - \frac{1}{2}$. The function denoted by $o(1)$ on the right hand side is actually $\ll_a (\log_2 q)^{-1/2}$.*

For any fixed $n \geq 0$, let
$$N(f;\delta,n) = N(f;\tfrac{1}{2}+\delta,\tfrac{n}{\lambda},\tfrac{n+1}{\lambda}), \quad N(\delta;n) = \sum_{f\in S_2(q)^*}^h N(f;\delta,n).$$

By integration by parts, taking $\psi$ to be any of the three functions indicated in Proposition 3, (note $\psi = \tfrac{1}{n^j}|\varphi^{(j)}|$, $j=0,1,2$) we have
$$\tfrac{1}{\lambda} Z_T(n) \leq \sum_{f\in S_2(q)^*}^h \sum_{\substack{\beta-1/2>\lambda^{-1} \\ n/\lambda \leq \gamma < (n+1)/\lambda}} \hat{\psi}(\lambda(\beta-\tfrac{1}{2}))$$
$$\leq N(\lambda^{-1};n)\hat{\psi}(1) + \lambda \int_{1/\lambda}^{1/2} N(\delta;n)\hat{\psi}'(\lambda\delta)d\delta.$$

Then we use Theorem 5 in the short range of $\delta$ ($\delta \leq \tfrac{\log_2 q}{\log q}$) and Theorem 4 for the remaining range, getting
$$\lambda \int_{1/\lambda}^{1/2} N(\delta;n)\hat{\psi}'(\lambda\delta)d\delta \leq \frac{a^2}{(1-a)^2} \int_{1/2}^{+\infty} (F(1,u)-F(a,u))\hat{\psi}'(u+\tfrac{1}{2})du + o_{a,\varphi}(\tfrac{1}{n^j}).$$

**Remark.** It can be verified using Lemma 2 below that with the choice of the test functions $\varphi_\varepsilon$ indicated in Step 5 above, the error terms encountered $o_{a,\varphi_\varepsilon}(\tfrac{1}{n^j})$ actually do not depend on $\varepsilon$.

**Lemma 1** *We have*
$$\frac{a^2}{(1-a)^2} \int_{1/2}^{+\infty} (F(1,u)-F(a,u))\hat{\psi}'(u+\tfrac{1}{2})du =$$
$$\frac{a^2}{(1-a)^2} \int_{-1}^{1} x\psi(x)e^{x/2}(K(1,x)-K(a,x))dx.$$

*Proof.* The integral can be transformed into
$$\int_{-1}^{1} x\psi(x)e^{x/2} \int_{1/2}^{+\infty} (F(1,u)-F(a,u))e^{xu}dudx$$
and then we use the following formulae to finish the computation:
$$\int_{1/2}^{+\infty} \frac{e^{-(\tfrac{2}{a}-x)u}}{u^2}du = 2E\Big(\tfrac{1}{2}(\tfrac{2}{a}-x)\Big) \tag{12}$$
$$\int_{1}^{+\infty} E(bu)\frac{e^{au}}{u}du = \frac{1}{a}\Big(E(b-a)-e^aE(b)\Big) \text{ for } b>a. \tag{13}$$

Of these, (12) is immediate from the definition (2), while (13) can be proved for instance by the following computation, using (2)
$$\int_{1}^{+\infty} E(bu)\frac{e^{au}}{u}du = \int_{1}^{+\infty} b\Big(\int_{1}^{+\infty} e^{-y}y^{-2}dy\Big)e^{au}du$$
$$= b\int_{b}^{+\infty} e^{-y}y^{-2} \int_{1}^{y/b} e^{au}dudy$$
$$= \frac{b}{a}\int_{b}^{+\infty} (e^{ay/b}-e^a)e^{-y}y^{-2}dy$$

and now the result is immediate.
□

From this the proposition follows.

### 3.4 Final estimates for the test function

The convergence of $\varphi_\varepsilon$ to $\varphi_0$ can be made slightly more precise.

**Lemma 2** *Let $h$ be any continuous function on $\mathbf{R}$. Then as $\varepsilon \to 0^+$, we have*

$$\int_\mathbf{R} |\varphi_\varepsilon(x)| h(x) dx \to \int_\mathbf{R} |\varphi_0(x)| h(x) dx$$

$$\int_\mathbf{R} |\varphi'_\varepsilon(x)| h(x) dx \to \int_\mathbf{R} |\varphi'_0(x)| h(x) dx$$

$$\int_\mathbf{R} |\varphi''_\varepsilon(x)| h(x) dx \to \int_\mathbf{R} |\varphi''_0(x)| h(x) dx$$
$$+ (|\varphi'_0|h)(-1^+) + (|\varphi'_0|h)(1^-) + (|\varphi'_0|h)(0^+) + (|\varphi'_0|h)(0^-)$$

*where the $\pm$ superscript indicate a limit from above or below.*

The proof is left to the reader.

**Lemma 3** *As $\varepsilon \to 0^+$, we have*

$$G_\psi(1) = 0.1535\ldots + o(1) \text{ for } \psi = |\varphi_\varepsilon|$$
$$G_\psi(1) = 0.3666\ldots + o(1) \text{ for } \psi = |\varphi'_\varepsilon|$$
$$G_\psi(1) = 0.3321\ldots + o(1) \text{ for } \psi = |\varphi''_\varepsilon|.$$

This follows from the previous lemma by direct computations.

Thus we are led to take $\psi = |\varphi_\varepsilon|$ for $n = 0, 1, 2$, and $\psi = |\varphi''_\varepsilon|/n^2$ for $n > 2$. From the previous steps (7) and (8), we get

$$\frac{1}{\lambda} Z'_T \leq 2 \frac{a^2}{(1-a)^2} \lim_{\varepsilon \to 0} \left\{ 3(G_{\varphi_\varepsilon}(1) - G_{\varphi_\varepsilon}(a)) + (\frac{\pi^2}{6} - \frac{5}{4})(G_{\varphi''_\varepsilon}(1) - G_{\varphi''_\varepsilon}(a)) \right\} + o_a(1)$$

which is Proposition 4.

## A  Proof of Theorem 4

Let $s = \sigma + it$ with $\frac{1}{2} \leq \sigma \leq 1$, and write $\sigma = \frac{1}{2} + \delta$. We have introduced in [KM2] the following "mollifier" $M(f, s)$, which is a smoothed partial sum of the Dirichlet series representing the inverse of $L(f, s)$: take $\Delta$ with $0 < \Delta < \frac{1}{2}$, $a$ with $0 < a < 1$, and let $M = q^\Delta$. Then we define first the continuous function $g_{M,a}$ by

$$g_{M,a}(x) = \begin{cases} 1 & \text{if } x \leq M^a \\ \dfrac{\log(x/M)}{\log(M^{a-1})} & \text{if } M^a \leq x \leq M \\ 0 & \text{if } x > M \end{cases} \quad (14)$$

and then we let (here $\varepsilon(n)$ is the trivial Dirichlet character modulo $q$)

$$x_m = \mu(m) m^{-\delta - it} \sum_{n \geq 1} \frac{\varepsilon(n) \mu(mn)^2}{n^{1+2\delta+2it}} g_{M,a}(mn)$$

and finally

$$M(f, s) = \sum_m \frac{x_m}{\sqrt{m}} \lambda_f(m).$$

Then we have, for $\sigma = \mathrm{Re}(s) > 1$

$$L(f,s) M(f,s) = 1 + \sum_{n > M^a} c_f(n) n^{-s} \qquad (15)$$

where the coefficients $c_f(n)$ satisfy

$$|c_f(n)| \leq \tau(n)^2.$$

Consider also the following arithmetic functions

$$\nu_\delta(k) = \frac{1}{k} \sum_{uv=k} \frac{\mu(u)}{u^{1+2\delta}} = \frac{1}{k} \prod_{p | k} (1 - p^{-(1+2\delta)})$$

$$\eta_t(n) = \sum_{uv=n} \left(\frac{u}{v}\right)^{it}$$

and the following coefficients, supported on squarefree integers $k$

$$y_k = \frac{\mu(k)}{k^{\delta+it}} \sum_{m,n} \frac{\mu(kmn)^2 \mu(m) \eta_t(m) n^{-it}}{(mn)^{1+2\delta+it}} g_{M,a}(kmn). \qquad (16)$$

The following lemma, when $\Delta < \frac{1}{4}$ and $a = \frac{1}{2}$, is contained in Section 3.4 of [KM2] (see also [Kow, 5.3]).

**Lemma 4** *For any $\Delta$ with $0 < \Delta < \frac{1}{2}$, any $0 < a < 1$, there exists an absolute constant $B > 0$ and a real number $\gamma = \gamma(\Delta) > 0$ such that we have*

$$\sideset{}{^h}\sum_{f \in S_2(q)^*} |L(f,s) M(f,s)|^2 \leq \zeta_q(1+2\delta) \sum_k \nu_\delta(k) |y_k|^2 + O(q^{-\gamma}(1+|t|)^B) \qquad (17)$$

*uniformly for $s$ with $\sigma - \frac{1}{2} \in [\frac{1}{\log q}, \frac{(\log_2 q)^2}{\log q}]$ and $t = \mathrm{Im}(s) \in \mathbf{R}$.*

As in the previous works ([Kow, Prop. 11]), we deduce from this lemma and from (15), by convexity, the

**Corollary 1** *With the same notations, for any $a'$ with $0 < a' < a$ and any $\sigma \geq \frac{1}{2}(\log q)^{-1}$, we have*

$$\sideset{}{^h}\sum_{f \in S_2(q)^*} |L(f,s) M(f,s) - 1|^2 \ll_{a,a',\Delta} (1+|t|)^B M^{-2a'(\sigma-\frac{1}{2})}.$$

Theorem 4 is derived from this following [Kow], [KM2] (see also below Appendix B).

We now sketch the proof of Lemma 4. The extension of the arguments of [KM2] to all $a$ is immediate; less clear is the extension to $\Delta < \frac{1}{2}$. This follows by an application of the methods of Iwaniec and Sarnak [IS]. Let $s = \sigma + it$ and

$$M_2 = \sum_{f \in S_2(q)^*}^{h} |L(f,s)M(f,s)|^2$$

the mollified second moment. As in [KM2, 3.3], we obtain an expression

$$\left(\frac{q}{4\pi^2}\right)^\delta H(\sigma - \tfrac{1}{2})M_2 = \sum_b \frac{1}{b} \sum_{n \geq 1} \sum_{m_1, m_2} \frac{\eta_t(n)}{\sqrt{m_1 m_2 n}} x_{bm_1} \bar{x}_{bm_2} U_s\left(\frac{4\pi^2 n}{q}\right) \Delta(m_1 m_2, n)$$

(see (25) for the definition of $\Delta(m,n)$), where $U_s$ is a certain function, similar to the function $V_s$ used below in Section C. The estimate of (25) for the remainder term $\mathcal{J}(m,n) = \Delta(m,n) - \delta(m,n)$ is sufficient to allow a mollifier of length $M = q^\Delta$ with $\Delta < \frac{1}{4}$. To go further, we use the explicit expression as a series of Kloosterman sums

$$\mathcal{J}(m,n) = -\frac{2\pi}{q} \sum_{r \geq 1} \frac{1}{r} S(m,n;qr) J_1\left(\frac{\sqrt{mn}}{qr}\right);$$

this is similar to what is done in [KM3], see also the Appendix to [Kow], and [IS]. This makes it possible to exploit the particular properties of the coefficients $\eta_t(n)$ involved in the sum so as, in effect, to detect cancellation when they are summed agains the Kloosterman sums, going beyond the Weil bound (which gives (25).)

Briefly, opening the Kloosterman sum and exchanging the order of summation of $n$ and $r$, we first study, for a fixed $a \bmod qr$, $(a, qr) = 1$, and a fixed $m < q$, the sum

$$\sum_{n \geq 1} \frac{\eta_t(n)}{\sqrt{n}} U_s\left(\frac{4\pi^2 n}{q}\right) J_1\left(\frac{\sqrt{mn}}{qr}\right) e\left(\frac{nx}{qr}\right).$$

The following summation formula is crucial.

**Lemma 5** *Let* $t : [0, +\infty[ \to \mathbf{C}$ *be a* $C^\infty$ *function, vanishing in a neighborhood of 0 and quickly decreasing at infinity. Let* $a, c \geq 1$ *be integers with* $(a, c) = 1$, *and* $d$ *the inverse of* $a$ *modulo* $c$. *For any* $t \neq 0$ *we have*

$$\begin{aligned}
\sum_{n \geq 1} \eta_t(n) e\left(\frac{an}{c}\right) t(n) &= \frac{\zeta(1+2it)}{c^{1+2it}} \int_0^{+\infty} t(x) x^{it} dt \\
&+ \frac{1}{c} \sum_{h \geq 1} \eta_t(h) e\left(-\frac{dh}{c}\right) \int_0^{+\infty} t(x) J_{2it}^+\left(\frac{4\pi\sqrt{hx}}{c}\right) dx \\
&+ \frac{1}{c} \sum_{h \geq 1} \eta_t(h) e\left(\frac{dh}{c}\right) \int_0^{+\infty} t(x) K_{2it}^+\left(\frac{4\pi\sqrt{hx}}{c}\right) dx
\end{aligned}$$

*where* $J^+$ *and* $K^+$ *are the following modified Bessel functions:*

$$J_{2it}^+(z) = -\frac{\pi}{\sin \pi it}(J_{2it}(z) - J_{-2it}(z)) \qquad (18)$$

$$K_{2it}^+(z) = 4\cos(\pi it) K_{2it}(z). \qquad (19)$$

This can be proved either by appealing to the modularity of the Eisenstein series with Fourier coefficients $\eta_t(n)$, or by classical abelian harmonic analysis (2-dimensionnal Poisson summation formula), see [D-I]. One can check that for $t \to 0$, this gives back the classical Voronoï formula used in [KM3] for instance.

After applying this transformation with

$$t(x) = \frac{1}{\sqrt{x}} U_s\left(\frac{4\pi^2 x}{q}\right) J_1\left(\frac{4\pi\sqrt{mx}}{qr}\right)$$

(more precisely, one has to multiply this by a test function $\xi$ which vanishes near 0 and is 1 for $x \geq 1$, as in [KM3], because of convergence difficulties due to the fact that the weight is 2), all the terms can be further estimated and shown to be small enough if $m < q$, which (since $m = m_1 m_2$ in the application to $M_2$) means that we can take a mollifier of length $M = q^\Delta$ for any $\Delta < \frac{1}{2}$, see [IS] for the full details in the case $t = 0$, [KM3] and [Kow, Ch. 6] for a similar problem, from which the techniques to do this can be adapted also.

**Remark.** The suspicious reader can simply take Lemma 4 in the state it is proved in [KM2], namely for $\Delta < \frac{1}{4}$ and compute that value of the explicit constant in Theorem 1 which is thereby obtained is 10.6 (for $a = 0.56$).

## B  Proof of Theorem 5

Theorem 5 is proved, like the other density theorem, by following the idea of Selberg [Sel, Lemma 14] (see [Kow, 5.2]); the next lemma provides the detector used to count zeros.

**Lemma 6** *Let $h$ be an holomorphic function in a neighborhood of the domain $\operatorname{Re}(s) \geq \sigma'$, $t_1 \leq \operatorname{Im}(s) \leq t_2$. Assume that in this region it satisfies*

$$h(s) = 1 + o\left(\exp\left(-\frac{\pi}{t_2 - t_1}\operatorname{Re}(s)\right)\right)$$

*uniformly for $\operatorname{Re}(s) \to +\infty$. Then, denoting the zeros of $h$ (with multiplicity) by $\rho = \beta + i\gamma$, we have*

$$2(t_2 - t_1) \sum_{\substack{\beta > \sigma' \\ t_1 < \gamma < t_2}} \sin\left(\pi \frac{\gamma - t_1}{t_2 - t_1}\right) \sinh\left(\pi \frac{\beta - \sigma'}{t_2 - t_1}\right) =$$

$$\int_{t_1}^{t_2} \sin\left(\pi \frac{t - t_1}{t_2 - t_1}\right) \log|h(\sigma' + it)| dt$$

$$+ \int_{\sigma'}^{+\infty} \sinh\left(\pi \frac{\beta - \sigma'}{t_2 - t_1}\right) \{\log|h(\beta + it_1)| + \log|h(\beta + it_2)|\} d\beta.$$

We apply this lemma to the functions $h(s) = L(f,s)M(f,s)$; the zeros of $L(f,s)$ are zeros of $h$ with at least the same multiplicity.

We will need two estimates given by the following propositions, which refine Lemma 4.

**Proposition 5** *With the notations and hypotheses of Lemma 4, there exists an absolute constant $B > 0$ such that we have for any $a' < a$ and any $\sigma \geq \frac{1}{2} + (\log q)^{-1}$ the estimate*

$$\sum_{f \in S_2(q)^*}^h \{|L(f,s)M(f,s)|^2 - 1\} \ll_{a,a',\Delta} (1 + |t|)^B M^{-a'(\sigma - \frac{1}{2})}.$$

**Proposition 6** *With notations as in Lemma 4, assume that $|t| \leq 1$ and*

$$\frac{1}{\log q} \leq \delta := \sigma - \tfrac{1}{2} \leq \frac{(\log_2 q)^2}{\log q}.$$

*Then we have as $q \to +\infty$*

$$\sideset{}{^h}\sum_{f \in S_2(q)^*} |L(f,s)M(f,s)|^2 \leq 1 + \frac{M^{-2a\delta} - M^{-2\delta}}{4\delta^2(1-a)^2(\log M)^2} + O_a\left(M^{-2a\delta}\frac{(\log_2 q)^3}{\log q}\right).$$

The first proposition is proved in Appendix C and the second one in Appendix D.

Now let $\sigma$ be such that

$$\frac{1}{\log q} \leq \delta := \sigma - \tfrac{1}{2} \leq \frac{\log_2 q}{\log q} \tag{20}$$

and $t_1$, $t_2$ such that $t_2 - t_1 = \lambda^{-1}$, where as in the statement of Theorem 5, we have $\lambda = a\Delta \log q$, and

$$-\frac{\log_2 q}{\log q} \leq t_1 < t_2 \leq \frac{\log_2 q}{\log q}.$$

We let

$$\sigma' = \sigma - \frac{1}{2\lambda}, \quad t_1' = t_1 - \frac{\mu}{\lambda}, \quad t_2' = t_2 + \frac{\mu}{\lambda}$$

where $\mu > 0$ is some parameter to be chosen later. For $f \in S_2(q)^*$ and any $\varepsilon > 0$, we have by (15)

$$L(f,s)M(f,s) - 1 = O_{\varepsilon,a}(M^{-a\mathrm{Re}(s)-1+\varepsilon})$$

uniformly for $\mathrm{Re}(s) \to +\infty$, and Lemma 6 can be applied in the region $\mathrm{Re}(s) \geq \sigma'$, $t_1' \leq \mathrm{Im}(s) \leq t_2'$, as soon as

$$a'\Delta \log q := \frac{\pi}{t_2' - t_1'} < a\Delta \log q$$

which means as soon as $2\mu + 1 > \pi$.

The zeros $\rho = \beta + i\gamma$ of $L(f,s)$, in the region $\beta \geq \sigma$, $t_1 \leq \gamma \leq t_2$, are among those of $h = LM$ in the above enlarged region. Moreover we see easily that for any such $\rho$ the inequality

$$1 \leq \frac{\lambda}{\pi \sin(\frac{\pi\mu}{2\mu+1})} 2(t_2' - t_1') \sinh\left(\pi \frac{\beta - \sigma'}{t_2' - t_1'}\right) \sin\left(\pi \frac{\gamma - t_1'}{t_2' - t_1'}\right)$$

holds: this is the detector. Summing over the forms $f \in S_2(q)^*$, we derive, by the inequality

$$\log |x| \leq \frac{|x|^2 - 1}{2},$$

that

$$\sideset{}{^h}\sum_{f \in S_2(q)^*} N(f; \sigma, t_1, t_2) \leq \frac{\lambda}{2\pi \sin(\frac{\pi\mu}{2\mu+1})} \Big\{ \int_{t_1'}^{t_2'} \sin\left(\pi \frac{\gamma - t_1'}{t_2' - t_1'}\right) \sideset{}{^h}\sum_{f \in S_2(q)^*} (|LM(f, \sigma' + it)|^2 - 1) dt$$

$$+ \int_{\sigma'}^{+\infty} \sinh\left(\pi \frac{\beta - \sigma'}{t_2' - t_1'}\right) \Big( \sideset{}{^h}\sum_{f \in S_2(q)^*} (|LM(f, \beta + it_1')|^2 - 1) + \sideset{}{^h}\sum_{f \in S_2(q)^*} (|LM(f, \beta + it_2')|^2 - 1) \Big) d\beta \Big\}$$

In the two last integrals, for $\beta$ in the range $\beta \geq \frac{1}{2} + \frac{(\log_2 q)^2}{\log q}$ we use Corollary 5. Those two last terms are then error terms, namely they are

$$\ll_\Delta \frac{1}{a - a'} e^{-\gamma(a-a')(\log_2 q)^2}, \text{ for some } \gamma := \gamma > 0 \tag{21}$$

Now to treat the main contributions (the first integral and the values of $\beta \leq \frac{1}{2} + \frac{(\log_2 q)^2}{\log q}$), we use the fine Lemma 6. We find that the two last integrals are actually convergent even if $2\mu + 1 = \pi$ (ie. $a' = a$), and are decreasing functions of $\mu$; so we henceforth replace $\mu$ by its smallest possible value, namely $2\mu + 1 = \pi$). For the first term we take $a' = a - \frac{1}{(\log_2 q)^{1/2}}$ so that (recall (20)) the error term in (21) is

$$O_{a,\Delta}(M^{-2a\Delta(\sigma - \frac{1}{2})} \frac{(\log_2 q)^3}{\log q}).$$

We get (recall that $c$ was defined in the statement of Theorem 5)

$$\sum_{f \in S_2(q)^*}^h N(f; \sigma, t_1, t_2) \leq \frac{a^2}{2c\lambda(1-a)^2} \Big\{ \frac{1}{(\sigma' - \frac{1}{2})^2} \int_{t_1'}^{t_2'} \sin\Big(\pi \frac{t - t_1'}{t_2' - t_1'}\Big) (e^{-2\lambda(\sigma' - \frac{1}{2})} - e^{-\frac{2}{a}\lambda(\sigma' - \frac{1}{2})}) dt$$

$$+ 2 \int_{\sigma'}^{+\infty} \sinh(\lambda(\beta - \sigma')) \frac{e^{-2\lambda(\beta - \frac{1}{2})} - e^{-\frac{2}{a}\lambda(\beta - \frac{1}{2})}}{(\beta - \frac{1}{2})^2} d\beta \Big\}.$$

$$+ O_a(M^{-2a\Delta(\sigma - \frac{1}{2})} \frac{(\log_2 q)^3}{\log q}).$$

Let $u = \lambda(\sigma - \frac{1}{2}) - \frac{1}{2}$. The first integral inside the brackets is equal to

$$(1 + O(a - a')) \frac{2\lambda}{u^2}(e^{-2u} - e^{-\frac{2u}{a}})$$

while the second integral is built from integrals of the form

$$I_a^\pm = e^{\mp u} \int_{\delta - 1/(2\lambda)}^{+\infty} \frac{e^{-(\frac{2}{a} \mp 1)\lambda\beta}}{\beta^2} d\beta$$

which are evaluated in terms of the function $E$:

$$I_a^\pm = \frac{\lambda e^{\mp u}}{u} E(u(\tfrac{2}{a} \mp 1)).$$

The computations yield explicitly

$$\sum_{f \in S_2(q)^*}^h N(f; \sigma, t_1, t_2) \leq \frac{a^2(1 + O_a((\log_2 q)^{-1/2}))}{(1-a)^2}(F(1, u) - F(a, u)) + O_a(\exp(-2u) \frac{(\log_2 q)^3}{\log q})$$

with

$$F(a, u) = \frac{1}{cu^2}\Big(e^{-2u/a} + ue^{-u} E(\frac{2-a}{a} u) - ue^u E(\frac{2+a}{a} u)\Big)$$

which concludes the proof of Theorem 5.

**Remark.** In [KM2], we had applied Selberg's lemma to the function $h(s) = 1 - (L(f,s)M(f,s) - 1)^2$, using the inequality

$$\log(1 + |x|^2) \leq |x|^2;$$

this procedure would have resulted in the loss of a factor 2 in the present case.

## C  Proof of Proposition 5

We begin with a lemma.

**Lemma 7** *With the same hypotheses as in Lemma 4, for any $a'$ with $0 < a' < a$, we have*

$$\sum_{f \in S_2(q)^*}^h (L(f,s)M(f,s) - 1) \ll_{a,a',\Delta} (1+|t|)^B M^{-a'(\sigma - \frac{1}{2})}$$

*for some absolute constant $B > 0$.*

*Proof.* First, for $\sigma > 2$ and $a' < a$, we have by (15)

$$\sum_{f \in S_2(q)^*}^h (L(f,s)M(f,s) - 1) \ll_{a'} M^{-a'(\sigma-1)}. \tag{22}$$

We now proceed as for the estimate of the second moment. Fix an integer $N \geq 3$ and a real polynomial $G$ which is such that

$$G(1-s) = G(s) \tag{23}$$
$$G(-1/2) = G(-3/2) = \ldots = G(-N+1/2) = 0. \tag{24}$$

By the usual contour-shifts and the functional equation of the $L$-function $L(f,s)$, we obtain for $f \in S_2(q)^*$ and $\frac{1}{2} < \mathrm{Re}(s) < 1$ the identity

$$G(s)\Gamma(s+\tfrac{1}{2})L(f,s) = \sum_{n \geq 1} \frac{\lambda_f(n)}{n^s} V_s\left(\frac{n}{\hat{q}}\right) + \varepsilon_f \hat{q}^{1-2s} \sum_{n \geq 1} \frac{\lambda_f(n)}{n^{1-s}} V_{1-s}\left(\frac{n}{\hat{q}}\right)$$

where we have denoted $\hat{q} = \frac{\sqrt{q}}{2\pi}$, and where

$$V_s(y) = \frac{1}{2i\pi} \int_{(2)} \Gamma(w+s+\tfrac{1}{2})G(w+s)y^{-w} \frac{dw}{w}.$$

The asymptotics of the function $V_s$ are easily determined by shifting the line of integration: we see that for some absolute constant $B > 0$, we have

$$V_s(y) = \Gamma(s+\tfrac{1}{2})G(s) + O(y^2(1+|t|)^B e^{-\pi|t|})$$

for $y \to 0$, and

$$V_s(y) = O(y^{-2}(1+|t|)^B e^{-\pi|t|})$$

for $y \to +\infty$.

These properties and an appeal to the Petersson formula, in the form of its corollary ([Du], [KM2] for instance)

$$\Delta(m,n) := \sum_{f \in S_2(q)^*}^h \lambda_f(m)\lambda_f(n) = \delta(m,n) + O_\varepsilon\left(\frac{(mn)^{1/2+\varepsilon}}{q^{3/2}}\right) \tag{25}$$

yield for the average of the values of the $L$-functions at $s$ the formula, for any $m < q$

$$\sum_{f \in S_2(q)^*}^h L(f,s)\lambda_f(m) = \frac{1}{m^s} + O_\varepsilon((1+|t|)^B m^{1/2} q^{-1/2+\varepsilon}).$$

Incorporating the mollifier, we deduce that there exists an absolute constant $B > 0$ such that for any $a$ with $0 < a < 1$ and any $\Delta$ with $0 < \Delta < \frac{1}{2}$, there exists $\gamma = \gamma(\Delta) > 0$ for which the estimate

$$\sum_{f \in S_2(q)^*}^h (L(f,s)M(f,s) - 1) \ll_{a,\Delta} (1+|t|)^B q^{-\gamma}$$

holds. The lemma follows from this and (22), again by an easy convexity argument.
□

The estimates for the first and second moments, Corollary 1 and Lemma 7, now obviously imply Proposition 5.

## D  Proof of Proposition 6

The starting point is still Lemma 4. We will estimate more precisely the sum

$$\sum_k \nu_\delta(k)|y_k|^2$$

to obtain the lemma. First note that from the definition (16), $y_k$ is supported on squarefree integers $k \leq M$. Recall also that one assume

$$\frac{1}{\log q} \leq \delta := \sigma - \tfrac{1}{2} \leq \frac{(\log_2 q)^2}{\log q}$$

**First case:**  Suppose first that $k \leq M^a$.

Under this hypothesis, as in [KM2, 3.4], we obtain the equality (note that $t$ has disappeared on the right-hand side) for $k$ squarefree

$$k^{\delta+it} y_k = \zeta_k(1+2\delta)^{-1} + \frac{1}{2i\pi} \int_\mathcal{C} \zeta_k(w+1+2\delta)^{-1} \frac{(M^a/k)^w (M^{(1-a)w} - 1)}{\log M^{1-a}} \frac{dw}{w^2}$$

where $\mathcal{C}$ is the curve defined by

$$\mathcal{C} := [-\frac{\kappa}{\log(U+2)} - iU, \ -\frac{\kappa}{\log(U+2)} + iU] \cup \{-\frac{\kappa}{\log(|t|+2)} + it, \ |t| \geq U\}$$

with $U = \exp((\log_2 q)^3)$, and $\kappa > 0$ an absolute constant such that $\zeta(1+s)$ admits no zeros on and to the right of $\mathcal{C}$. Notice that $\kappa/\log(U+2) > 4\delta$. The classical estimates of Hadamard and de la Vallée Poussin (see [Tit, Ch. 3]) for $\zeta$ and $\zeta^{-1}$ can be written in the form

$$|\zeta(1+s)| \leq K \log(|\text{Im}(s)| + 1), \quad |\zeta(1+s)^{-1}| \leq K \log(|\text{Im}(s)| + 1)$$

for some absolute constant $K > 0$, for all $s \in \mathcal{C}$.

Let
$$\omega_s(k) = \prod_{p|k} (1 - p^{-s})^{-1}.$$

For $k \leq M^a$, $k$ squarefree, the above formula thus yields
$$\left| y_k - \frac{\omega_{1+2\delta}(k)}{k^{\delta+it}} \zeta(1+2\delta)^{-1} \right| \ll \frac{1}{k^\delta \log M^{(1-a)}} \omega_{3/4}(k) \left( \left( \frac{k}{M^a} \right)^{\kappa/\log U} + \frac{\log U}{U} \right)$$

It follows that
$$\zeta_q(1+2\delta) \sum_{k \leq M^a} \nu_\delta(k) |y_k|^2$$
$$= \zeta(1+2\delta)^{-1} \sum_{k \leq M^a} \frac{\mu(k)^2 \omega_{1+2\delta}(k)}{k^{1+2\delta}} + O_a\left( M^{-2a\delta} \left( \frac{\log U}{\kappa \log M} + \frac{\log U}{U \log M} \right) \right). \qquad (26)$$

By our hypothesis on $\delta$ and our choice of $U$ the last error term is
$$O_a(M^{-2a\delta} \frac{(\log_2 q)^3}{\log q}).$$

Next comes the case $M^a < k < M$. In this case we use the integral representation
$$k^{\delta+it} y_k = \frac{1}{2i\pi} \int_{(2)} \zeta_k(w+1+2\delta)^{-1} \frac{(M/k)^s}{\log M^{1-a}} \frac{dw}{w^2}$$

which, after shifting the contour to $\mathcal{C}$, yields
$$k^{\delta+it} y_k = R_k + \frac{1}{2i\pi} \int_{\mathcal{C}} \zeta_k(w+1+2\delta)^{-1} \frac{(M/k)^s}{\log M^{1-a}} \frac{dw}{w^2}$$

where $R_k$ is the residue at $w = 0$, where the integrand has a double pole. Hence we compute
$$R_k = \zeta(1+2\delta)^{-1} \omega_{1+2\delta}(k) \frac{\log(M/k)}{\log M^{1-a}} + \frac{1}{\log M^{1-a}} \frac{d}{dw} \zeta_k(w+1+2\delta)^{-1} \Big|_{w=0}$$
$$= \zeta(1+2\delta)^{-1} \frac{\omega_{1+2\delta}(k)}{\log M^{1-a}} \left( \log \frac{M}{k} - \frac{\zeta'}{\zeta}(1+2\delta) - \sum_{p|k} \frac{\log p}{p^{1+2\delta}} \omega_{1+2\delta}(p) \right).$$

The third term in the inner sum is $\ll \omega_{3/4}(k)$ hence gives a negligible contribution because of our hypothesis on $\delta$. We get the estimate
$$\left| y_k - \frac{\zeta(1+2\delta)^{-1}}{k^{\delta+it}} \frac{\omega_{1+2\delta}(k)}{\log M^{1-a}} (\log M/k - \frac{\zeta'}{\zeta}(1+2\delta)) \right|$$
$$\ll \frac{\omega_{3/4}(k)}{\log M^{(1-a)}} \left( \left( \frac{k}{M^a} \right)^{\kappa/\log U} + \frac{\log U}{U} + \frac{\omega_{3/4}(k)}{\zeta(1+2\delta)} \right)$$

which gives therefore

$$\zeta_q(1+2\delta)\sum_{M^a<k\leq M}\nu_\delta(k)|y_k|^2 =$$
$$\frac{\zeta(1+2\delta)^{-1}}{(\log M^{1-a})^2}\sum_k\frac{\mu(k)^2\omega_{1+2\delta}(k)}{k^{1+2\delta}}\Big(\log\frac{M}{k}-\frac{\zeta'}{\zeta}(1+2\delta)\Big)^2 \qquad (27)$$
$$+O_a\Big(M^{-2a\delta}\big(\frac{\log U}{\kappa\log M}+\frac{\log U}{U\log M}\big)\Big)$$

and once again the error term is

$$O_a(M^{-2a\delta}\frac{(\log_2 q)^3}{\log q}).$$

Therefore, we have to estimate the sum $S$ defined by

$$S = \zeta(1+2\delta)^{-1}\Big\{\sum_{k\leq M^a}\frac{\mu(k)^2\omega_{1+2\delta}(k)}{k^{1+2\delta}}$$
$$+\frac{1}{(\log M^{1-a})^2}\sum_{M^a<k\leq M}\frac{\mu(k)^2\omega_{1+2\delta}(k)}{k^{1+2\delta}}\Big(\log\frac{M}{k}-\frac{\zeta'}{\zeta}(1+2\delta)\Big)^2\Big\}. \qquad (28)$$

For this we first note the identity

$$\frac{2}{\log M^{1-a}}\frac{1}{2i\pi}\int_{(2)}\frac{M^{aw}}{k^w}\Big(\frac{M^{(1-a)w}-1}{w\log M^{1-a}}-1\Big)\frac{dw}{w^2} = \begin{cases} 1, & \text{if } k<M^a \\ \frac{\log^2 M/k}{(\log M^{1-a})^2} & \text{if } M^a<k\leq M \\ 0 & \text{if } k>M \end{cases}$$

We consider first

$$S_1 = \zeta(1+2\delta)^{-1}\Big\{\sum_{k\leq M^a}\frac{\mu(k)^2\omega_{1+2\delta}(k)}{k^{1+2\delta}}+\frac{1}{(\log M^{1-a})^2}\sum_{M^a<k\leq M}\frac{\mu(k)^2\omega_{1+2\delta}(k)}{k^{1+2\delta}}\Big(\log\frac{M}{k}\Big)^2\Big\}.$$

We have, by the above identity,

$$S_1 = \frac{2\zeta(1+2\delta)^{-1}}{\log M^{1-a}}\frac{1}{2i\pi}\int_{(2)}\zeta(1+2\delta+w)M^{aw}\Big(\frac{M^{(1-a)w}-1}{w\log M^{1-a}}-1\Big)\frac{dw}{w^2}$$

and we shift the contour to $\text{Re}(s)=-\frac{1}{2}$, passing by two poles at $s=0$ and $s=-2\delta$, which gives (using moreover $\zeta(1+2\delta)^{-1}=2\delta+O(\delta^2)$)

$$S_1 = 1+\frac{1}{\delta\log M^{1-a}}\Big(\frac{M^{-2a\delta}-M^{-2\delta}}{2\delta\log M^{1-a}}-M^{-2a\delta}\Big)+O\Big(\delta M^{-2a\delta}\Big). \qquad (29)$$

Similarly we consider the sum

$$S_2 = -2\frac{(\zeta'\zeta^{-2})(1+2\delta)}{(\log M^{1-a})^2}\sum_{M^a<k\leq M}\frac{\mu(k)^2\omega_{1+2\delta}(k)}{k^{1+2\delta}}\log\frac{M}{k}$$

We have the identity

$$\frac{1}{2i\pi}\int_{(2)}\frac{M^{aw}}{k^w}\Big(\frac{M^{(1-a)w}-1}{\log M^{1-a}}\Big)\frac{dw}{w^2} = \begin{cases} 1, & \text{if } k < M^a \\ \frac{\log M/k}{\log M^{1-a}} & \text{if } k > M^a \\ 0 & \text{if } k > M \end{cases}$$

and we find

$$\begin{aligned}S_2 &= -2\frac{(\zeta'\zeta^{-2})(1+2\delta)}{\log M^{1-a}}\frac{1}{2i\pi}\int_{(2)}\zeta(1+2\delta+w)M^{aw}\frac{M^{(1-a)w}-1}{\log M^{1-a}}\frac{dw}{w^2}\\ &\quad +2\frac{(\zeta'\zeta^{-2})(1+2\delta)}{\log M^{1-a}}\sum_{k<M^a}\frac{\mu(k)^2\omega_{1+2\delta}(k)}{k^{1+2\delta}}\end{aligned}$$

The first term on the right equals

$$-2\frac{(\zeta'\zeta^{-2})(1+2\delta)}{\log M^{1-a}}\Big[\zeta(1+2\delta)+\frac{M^{-2\delta}-M^{-2a\delta}}{4\delta^2\log M^{1-a}}\Big]$$

To treat the second term above, we use the following identity : let $\eta = 1/\log^{100} q$, then for any $M' > 1$, not an integer, we have

$$\frac{1}{2i\pi}\int_{(2)}\frac{M'^w}{k^w}\frac{\eta^{-1}dw}{(w+\eta^{-1})w} = \begin{cases} 1-(\frac{k}{M'})^{1/\eta}, & \text{if } k < M' \\ 0 & \text{if } k > M' \end{cases}$$

from this we infer

$$\sum_{k<M'}\frac{\mu(k)^2\omega_{1+2\delta}(k)}{k^{1+2\delta}} = \zeta(1+2\delta)-\frac{M'^{-2\delta}}{2\delta}+O(\eta M'^{-2\delta})+O(\eta^{-1}M'^{-1/2}) \qquad (30)$$

So we obtain

$$S_2 = -2\frac{(\zeta'\zeta^{-2})(1+2\delta)}{\log M^{1-a}}\Big[\frac{M^{-2\delta}-M^{-2a\delta}}{4\delta^2\log M^{1-a}}+\frac{M^{-2a\delta}}{2\delta}\Big]+O(\log^{-100} q M^{-2a\delta}) \qquad (31)$$

The last sum is

$$S_3 = \frac{(\zeta'^2\zeta^{-3})(1+2\delta)}{(\log M^{1-a})^2}\sum_{M^a<k\leq M}\frac{\mu(k)^2\omega_{1+2\delta}(k)}{k^{1+2\delta}}$$

which we find to be

$$S_3 = \frac{(\zeta'^2\zeta^{-3})(1+2\delta)(M^{-2a\delta}-M^{-2\delta})}{2\delta(\log M^{1-a})^2}+O\Big(\frac{M^{-2a\delta}}{\log^{100} q}\Big).$$

From the definition of $S$, with (29), (31) and this last estimate, we obtain

$$S = 1+\frac{M^{-2a\delta}-M^{-2\delta}}{4\delta^2(1-a)^2(\log M)^2}+O\Big(\delta M^{-2a\delta}\Big). \qquad (32)$$

Proposition 6 is now proved (see (27), (28)).

# References

[Br1] Brumer, A.: The rank of $J_0(N)$, Astérisque 228, SMF (1995), 41–68.

[D-I] Deshouillers, J.M. and Iwaniec, H.: Kloosterman sums and Fourier coefficients of cusp forms, Invent. Math. 70 (1982), 219–288.

[Du] Duke, W.: The critical order of vanishing of automorphic $L$-functions with high level, Invent. Math. 119, (1995), 165–174.

[IS] Iwaniec, H. and Sarnak, P.: The non-vanishing of central values of automorphic $L$-functions and Siegel-Landau zeros, in preparation.

[KM1] Kowalski, E. and Michel, P.: Sur le rang de $J_0(q)$, Préprint de l'Université d'Orsay, (1997).

[KM2] Kowalski, E. and Michel, P.: Sur les zéros des fonctions $L$-automorphes de grand niveau, Préprint de l'Université d'Orsay, (1997).

[KM3] Kowalski, E. and Michel, P.: A lower bound for the rank of $J_0(q)$, Préprint de l'Université d'Orsay, (1997).

[Kow] Kowalski, E.: The rank of the Jacobian of modular curves: analytic methods, Ph.D. Thesis, Rutgers University, 1998.

[MSD] Mazur, B. and Swinnerton-Dyer, P.: Arithmetic of Weil curves, Invent. Math., 25 (1974), 1–61.

[Tit] Titchmarsh, E.C.:The theory of the Riemann zeta function, 2nd Edition (revised by D.R. Heath-Brown), Oxford University Press, 1986.

[Sel] Selberg, A.: Contributions to the theory of Dirichlet's $L$-functions, Skr. Norske Vid. Akad. Oslo. I. (1946), 1–62, or Collected Papers, vol. 1, Springer Verlag, Berlin, (1989), 281–340.

P. Michel: michel@darboux.math.univ-montp2.fr
Université Montpellier II cc 051, 34095 Montpellier Cedex 05, FRANCE.

E. Kowalski: ekowalsk@math.princeton.edu
Dept. of Math. - Fine Hall, Princeton University, Princeton, NJ 08544-1000, USA.